# A New Initial Approximation Bound in the Durand-Kerner Algorithm for Finding Polynomial Zeros

Bandung Arry Sanjoyo[1*], Mahmud Yunus[1] and Nurul Hidayat[1]

[1] Departements of Mathematics, Institut Teknologi Sepuluh Nopember, Surabaya, Indonesia
* Corresponding author:bandung@matematika.its.ac.id

**Abstract.** The Durand-Kerner algorithm is a widely used iterative technique for simultaneously finding all the roots of a polynomial. However, its convergence heavily depends on the choice of initial approximations. This paper introduces two novel approaches for determining the initial values: New bound 1 and the lambda maximal bound, aimed at improving the stability and convergence speed of the algorithm. Theoretical analysis and numerical experiments were conducted to evaluate the effectiveness of these bounds. The lambda maximal bound consistently ensures that all the roots lie within the complex circle, leading to faster and more stable convergence. Comparative results demonstrate that while New bound 1 guarantees convergence, but it yields excessively large radii.

**Keywords:** Polynomial Zeros, Durand-Kerner, Initial Approximation, Root-Finding Algorithm.

## 1 Introduction

The problem of finding all the roots of a polynomial of degree $n$ is fundamental in mathematics, computer science, and various engineering applications [1], [2]. One of the widely used methods for solving this problem is the Durand-Kerner algorithm, also known as the Weierstrass method, which iteratively approximates all roots simultaneously. Although the Durand-Kerner method is popular due to its simplicity and parallelizable structure, it may fail to converge or yield inaccurate results, particularly for high-degree polynomials [3], [4], [5], [6].

Numerous studies have been conducted to improve the convergence speed and accuracy of the root-finding algorithms. One of the key strategies involves refining the selection of initial approximation to ensure stability and convergence [3], [7], [8], [9], [10], [11], [12]. This paper proposes a new approach for determining the initial approximation of the radius in the complex plane, aiming to enhance the performance of the Durand-Kerner method. The proposed bounds were evaluated through theoretical analysis and numerical experiments.

## 2 State of the Art and Related Work

### 2.1 Durand-Kerner Algorithm

A general monic polynomial zeros with real coefficients can be expressed as

$$p(x) = x^n + a_2 x^{n-1} + a_3 x^{n-2} + \cdots + a_n x + a_{n+1} = 0 \tag{1}$$

where $a_i \in \mathbb{R}$ [13], [14]. The Equation (1) has $n$ complex roots, or $x_i \in \mathbb{C}$. When the polynomial equation is expressed in its linear factorized form, the roots can be directly identified from the corresponding linear factors. Consequently, we get the following formula in Equation (2).

$$x_i^{(k+1)} = x^k - \frac{p(x)}{\prod_{j \neq i}^n (x - x_j)} \tag{2}$$

The convergence of this method is highly sensitive to the choice of initial approximations $x_i^0$ [14], [15], [16].

The Durand-Kerner algorithm is an iterative method based on Equation (2) for simultaneously computing all complex roots $x_i \in \mathbb{C}$, where $i = 1, 2, \ldots, n$. The method begins with initial approximations $x_i^0$ in the complex plane that are sufficiently close to the actual roots. The roots are then updated iteratively using Equation (2), so that the values $x_i^{k+1}$ that lead to and are very close to $x_i$. Globally, the steps of root finding are written in Algorithm 1.

---

Input: $coeffs = [1, a_2, a_3, \ldots, a_{n+1}]$ which is the coefficient of $p(x)$.
Output: $x = [x_1, x_2, \ldots, x_n]$ which is the roots of $p(x) = 0$.

Algorithm:
1. Set initial value $x_i^0$, i=1, 2, ..., n.
2. Compute next $x_i^{(k+1)} = x^k - \frac{p(x^k)}{\prod_{j \neq i}^n (x^k - x_j)}$, i=1, 2, ..., n.
3. Repeat step 2. until $x_i^{(k+1)}$ closed to $x_i^k$.

---

*Algorithm 1. Durand-Kerner algorithm.*

Giving the initial value $x_i^0$, it should be inside or on a circle of complex fields and close enough to the roots $x_i$ [17]. The algorithm converges when either the difference between successive approximations $|x_i^{k+1} - x_i^k|$ is less than a tolerance $\epsilon_1$, or the polynomial value $p(x_i^{k+1})$ is less than $\epsilon_2$, or the number of iterations exceeds a maximum threshold. The computational complexity of the Durand-Kerner method is $O(kn^2)$ where $k$ is the number of iterations and $n$ is the degree of the polynomial. The $O(kn^2)$ work comes from steps 2 and 3.

## 2.2 Setting the Radius of the Complex Plane

The choice of initial approximations significantly affects the convergence behavior of the Durand-Kerner method. Several studies have proposed bounds for the radius of the complex plane within which all roots are guaranteed to lie. Kjellberg [17] recommended that the initial approximations be greater than the absolute value of all polynomial roots [18], [19], [20], [21]. These bounds serve as a guide for selecting the initial approximations.

1. **Cauchy's Bound**

Cauchy proposed a bound where all roots lie within a circle of radius in Equation (3).

$$r = 1 + \max_{1 \leq i \leq n} | a_i | \qquad (3)$$

The procedure for determining the radius described in Equation (3) is outlined in Algorithm 2.

> Input: $coeffs = [1, a_2, a_3, \ldots, a_{n+1}]$ which is the coefficient of $p(x)$.
> Output: $r$, where $r$ denotes the radius of the complex plane that contains the initial approximations $x_i^0$.
>
> Algorithm:
> ```
> function r=cauchyBound(coeffs)
>     r = (1+max(abs(coeffs)));
> end
> ```

*Algorithm 2. Algorithm for finding the radius of Cauchy;s bound.*

The *cauchyBound* algorithm requires constant extra space and has a computational complexity of $O(n)$ floating-point operations (flops).

2. **Lagrange's Bound**

Lagrange refined the bound on the roots by considering:

$$r = 1 + \max_{1 \leq i \leq n} | a_i |^{\frac{1}{i}} \qquad (4)$$

For the cases $|a_i| \geq 1$, the value of Cauchy's bound radius is greater than the value of the Lagrange's bound radius. In many cases $|a_i| < 1$, the value of Cauchy's bound radius is smaller than the value of the Lagrange's bound radius. Algorithm 3 provides the steps for computing the radius as defined in Equation (4).

```
Input: coeffs = [1, a_2, a_3, ..., a_{n+1}] which is the coefficient of p(x).
Output: r, where r denotes the radius of the complex plane that contains the
        initial approximations x_i^0.

Algoritma:
function r= lagrangeBound (coeffs)
    maksLagrange=0;
    for k=2:n
        M=abs(coeffs(k))^(1/k);
        if maksLagrange < M
            maksLagrange = M;
        end
    end
    r = (1 + maksLagrange);
end
```

*Algorithm 3. Algorithm for finding the radius of Lagrange's bound*

The *lagrangeBound* algorithm requires constant extra space and has a computational cost of $O(n^2)$ flops. This quadratic workload does not increase the overall complexity of Algorithm 1.

3. **Aberth's Bound**

Aberth introduced a bound based on a modified polynomial transformation, which also ensures all roots are enclosed within a certain radius $r$ as in Equation (5) [13].

$$r = \frac{a_2}{n} + r_0 \quad (5)$$

where $r_0$ is a positive integer $w$ such that the value of $s(w)$, as defined in Equation (7), is positive. The function $s(w)$ is obtained from Equation (6) by replacing the plus sign on the coefficient $c_i$ with minus sign, as in Equation (7).

$$p\left(w - \frac{a_2}{n}\right) = w^n + c_2 w^{n-2} + c_3 w^{n-3} + \cdots + c_n \quad (6)$$

$$s(w) = w^n - c_1 w^{n-2} - c_2 w^{n-3} - \cdots - c_0 \quad (7)$$

The algorithm for computing the radius based on Equations (5), (6), and (7) is presented in Algorithm 4.

Input: $coeffs = [1, a_2, a_3, ..., a_{n+1}]$ which is the coefficient of $p(x)$.
Output: $r$, where $r$ denotes the radius of the complex plane that contains the initial approximations $x_i^0$.

Algorithm:
```
function r= aberthBound(coeffs)
  n1=-coeffs(2)/n; r=1;
  sw=c(n+1); //syms w;
  for i=1:n
      sw=sw+ coeffs(i)*(w+n1)^(n-i+1);
  end
  sw=expand(sw);
  cc = sym2poly(sw);
  cc(2:end)=-abs(cc(2:end));
  i=0; Maks=max(abs(coeffs),n);
  while (i<=Maks)
     if (polyval(cc,i) > 0)
        r=i;
        i=Maks+1;
     end
     i=i+1;
  end
end
```

*Algorithm 4. Algorithm for finding the radius of Aberth's bound*

The *aberthBound* algorithm requires constant extra space and has a computational cost of $O(n^2)$ flops.

Each of these bounds has its own computational characteristics and accuracy. However, they may not always provide the tightest possible radius, especially for polynomials with specific coefficient structures. This motivates the development of new bounds, such as those proposed in this paper, to improve convergence reliability and efficiency.

## 3  New Bound Proposed

To determine a more effective radius of the complex plane that encloses all the roots of a polynomial, both theoretical analysis and numerical experiments were conducted. The experimental observations were based on the graphical analysis of $p(x)$, using function plotting software and specific polynomial structures.

**Theorem 1**
For a monic polynomial $p(x) = x^n + a_{n+1}$, where $a_i = 0$ for $i = 2, 3, ..., n$, all the roots lie within a complex circle of radius:

$$r = \left|(a_{n+1})^{\frac{1}{n}}\right| \tag{8}$$

**Theorem 2**
For a monic polynomial in Equation (1), if $a_s \neq 0$ and $a_i = 0$ for $2 \leq i \neq s \leq (n+1)$, then the roots lie within a circle of radius all the roots lie within a complex circle of radius:

$$r = \left|(a_k)^{\frac{1}{k-1}}\right| \tag{9}$$

From the general form of the polynomial, and by applying the triangle inequality to the modulus of both sides of the polynomial equation, an inequality (11) is derived.

$$x^n = -(a_2 x^{n-1} + a_3 x^{n-2} + \cdots + a_n x + a_{n+1}) \tag{10}$$

Letting $|x| = r$, and dividing both sides of Equation (10) by $r^{n-1}$, we obtain inequality (11).

$$r \leq |a_2| + \sum_{k=3}^{n+1} |a_k| \frac{1}{r^{k-2}} \tag{11}$$

This inequality still contains the variable $r$ on the right-hand side. To resolve this, a comparison function is constructed as follows:

$$f(r, a) = |a_2| + \sum_{k=3}^{n+1} |a_k| \frac{1}{r^{k-2}} \tag{12}$$

and

$$g(r, a) = |a_2| + \sum_{k=3}^{n+1} |a_k|^{k-1} \tag{13}$$

Equation (12) and (13) indicate that the functions $f(r, a)$ and $g(r, a)$ are primarily influenced by the coefficient polynomial $a_2$. A comparison between $f(r, a)$ and $g(r, a)$ is illustrated in Fig. 1. In Fig. 1 (a) show that for $|a_2| < 1$ and $r > 2.5$, both functions lie below the modulus of the maximum root. In Fig. 1 (b) and (c), for $a_2 \geq 1$, the function $g(r, a)$ provides a tighter bound. Therefore, the following expression is proposed as a *New Bound 1*:

$$r = \sum_{k=2}^{n} \left|(a_k)^{\frac{1}{k-1}}\right| \tag{14}$$

It can be observed that the formulas in Theorems 1 and 2 are special cases of Equation (14). A comparison of the complex plane radius bounds given in Equations (3), (4), (5), and (14) is presented in the next section..

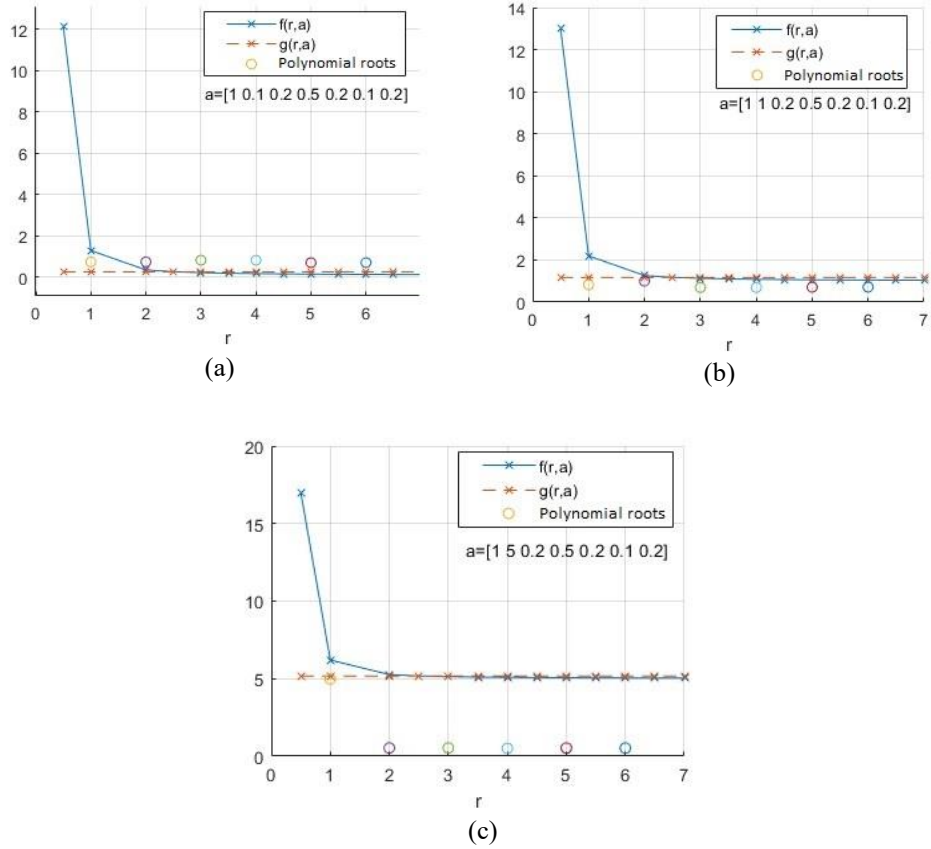

Fig. 1. A comparison between $f(r,a)$ and $g(r,a)$

**Definition 1: Companion Matrix**
Given a monic polynomial:

$$p(\lambda) = \lambda^n + a_2\lambda^{n-1} + a_3\lambda^{n-2} + \cdots + +a_n\lambda + a_{n+1} \qquad (15)$$

Its companion matrix $C \in \mathbb{C}^{n \times n}$ is defined as [22]

$$C = \begin{bmatrix} -a_{n-1} & -a_{n-2} & -a_{n-3} & \cdots & -a_0 \\ 1 & 0 & 0 & \cdots & 0 \\ 0 & 1 & 0 & \cdots & 0 \\ \vdots & \vdots & \ddots & 0 & 0 \\ 0 & 0 & 0 & 1 & 0 \end{bmatrix} = \begin{bmatrix} b & a_0 \\ I & 0 \end{bmatrix} \qquad (16)$$

where $I$ is identity matrix size $(n-1) \times (n-1)$ and $b = [a_{n-1} \quad a_{n-2} \quad \cdots \quad a_1]$. Let $|\lambda_1|$ be the dominant eigenvalue of $C$, i.e., the eigenvalue with the largest magnitude.

**Theorem 3: Lambda Maximal Bound**

If $\lambda_1$ is dominant eigenvalue of the companion matrix $C$, then all the roots of the polynomial lie within a complex circle of radius:

$$r = |\lambda_1| \tag{17}$$

Since $\lambda_1$ is the dominant eigenvalue, all the roots of the polynomial $p(x)$ lie within a complex circle of radius $r = \lambda_1$. This bound is the tightest among those considered, including Cauchy, Lagrange, Aberth, and the proposed summation bound. There for, it is recommended as the initial approximation of the radius for the Durand-Kerner method.

To construct the algorithm for determining the radius using the dominant eigenvalue, the power method is employed, defined as follows:

$$\lambda_1 = \frac{x_k^T x_{k+1}}{x_k^T x_k} \tag{18}$$

where $x_0 \in \mathbb{R}^{n \times 1}$, and $x_k = C^k x_0$.

The process of determining the radius $r$, as defined in Equation (14), is expressed in Algorithm 5. The *lambdaMaximal* algorithm requires non-constant extra space and has a computational complexity of $O(n^2)$ flops.

Input: $coeffs = [1, a_2, a_3, \ldots, a_{n+1}]$ which is the coefficient of $p(x)$.
Ouput: $r$, where $r$ denotes the radius of the complex plane that contains the initial approximations $x_i^0$.

Algorithm:
```
function r= lambdaMaximal(coeffs)
  C=[-coeffs(2:end); [eye(n-1) zeros(n-1,1)]];
  x0=2*rand(n,1);
  for i=1:20
     x1=C*x0; x=x0;
     x0=x1;
  end
  r = abs((x'*x1)/(x'*x));
end
```

*Algorithm 5. Algorithm for finding the radius of lambda maximal bound*

The *lambdaMaximal* algorithm requires an extra space is not constant, and the amount of work done is $O(n^2)$ flops.

## 4 Numerical Results and Discussion

### 4.1 Comparison of Radius Bounds

In numerical experiments, 50 monic polynomials of degrees ranging from 3 to 50 were generated with random coefficients in the range [−15,15]. For each polynomial, the radius bounds were computed using four methods: Cauchy, Lagrange, Aberth, and the proposed New Bound 1 and Lambda Maximal Bound.

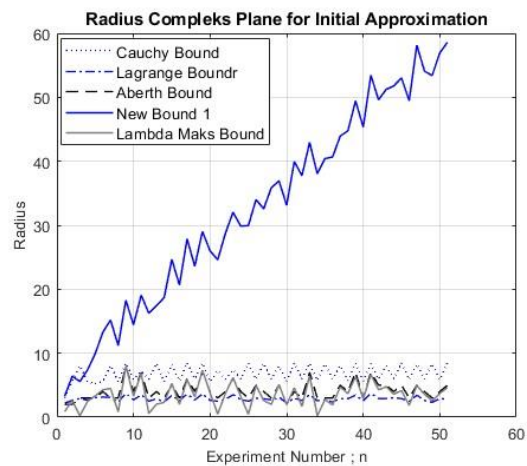

*Fig. 2. Comparison of the radius of the complex plane of the Cauchy, Lagrange, and Aberth boundaries with the proposed boundaries of New Bound 1 and Lambda dominant.*

Fig. 2 shows the radius comparison. New Bound 1 consistently produced the largest radii among the methods. While this ensures that all roots are enclosed, it may lead to slower convergence due to the initial guesses being far from the actual roots. The lambda maximal bound, derived from the dominant eigenvalue of the companion matrix, yielded the smallest radius that still enclosed all the roots. This resulted in faster convergence and fewer iterations. The Cauchy, Lagrange, and Aberth bounds showed moderate performance but occasionally failed to enclose all roots, as seen in the case of the 3rd degree polynomial in *Fig. 4*. *Fig. 4* shows the process of convergence of $x_i^k$ for some polynomial equations. The root-finding process for $x_i$, using the New Bound 1 limit, also ensures convergence.

The convergence behavior is displayed in Fig. 3 and Fig. 4. The convergence plots illustrate that the lambda maximal bound leads to rapid convergence with fewer iterations. New bound 1, despite its large radius, also converges but requires more iterations, making it less efficient. The lambda maximal bound strikes a balance between tightness and reliability, making it the most practical choice for initializing the Durand-Kerner algorithm.

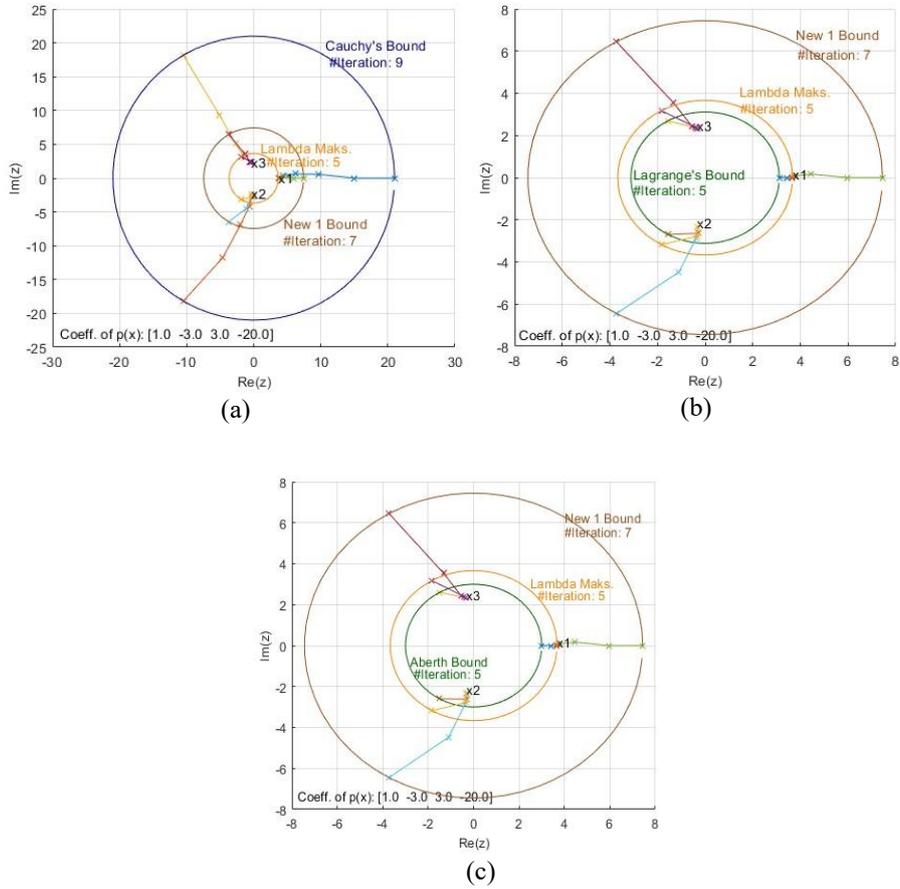

(a) (b)

(c)

*Fig. 3. Comparison of the convergence of the roots using proposed bound and classic bound.*

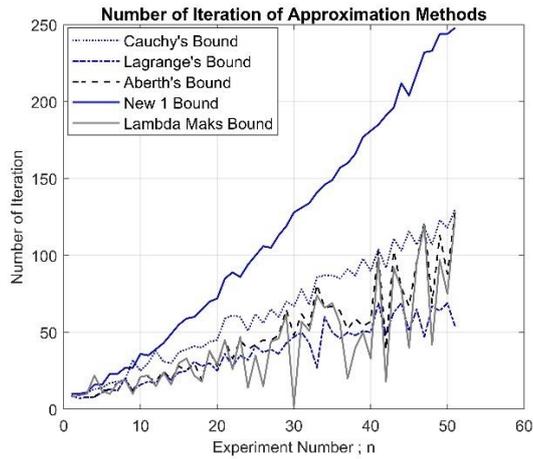

*Fig. 4. Graphical of the convergence behavior of $x_i^k$ under the proposed method.*

## 4.2 Accuracy and Convergence Behavior of the Durand-Kerner Algorithm

This section analyzes how initial radius bounds of new bounds influence the precision of the root estimates and the convergence. We evaluate the impact of New Bound 1, Lambda Maximal Bound ($\lambda_{max}$) in comparison with Matlab's built-in roots() function through a series of numerical experiments conducted on Wilkinson polynomials, randomly generated polynomials, and polynomials with clustered roots.

**Analysis on Diverse Polynomial Structures**

In the numerical experiments involving polynomials with well-separated roots, we first employed the Wilkinson polynomial defined as $\prod_{i=1}^{n}(x+i)$. The accuracy of the Durand-Kerner algorithm using initial approximations derived from New Bound 1 and the Lambda Maximal Bound ($\lambda_{max}$) is illustrated in the plots presented in Table 1. For polynomials of degree up to $n = 20$, both bounds yielded high accuracy, with a mean error of $e_{mean} = 1.193 \times 10 - 3$, which is significantly better than the mean error produced by MATLAB's built-in *roots*() function, recorded at $2.052 \times 10 - 2$. For higher-degree polynomials ($n \geq 40$), the mean error increased to the order of one decimal digit yet still outperformed the accuracy of the *roots*() function.

The convergence experiments were conducted using the Wilkinson polynomial defined as $\prod_{i=1}^{n}(x+i)$. The results, illustrated in the plots presented in Table 2, evaluate the convergence behavior of the Durand-Kerner algorithm when initialized with New Bound 1 and the $\lambda_{max}$ Bound. For polynomial degrees up to $n = 40$, both bounds successfully guided the algorithm to converge to the true roots, comparable to the performance of MATLAB's built-in roots() function. Furthermore, even for higher-degree polynomials up to $n = 140$, the algorithm maintained convergence, demonstrating the robustness of the proposed initial bounds in supporting stable root-finding across a wide range of polynomial complexities.

The Wilkinson polynomials are known to be highly sensitive to even little perturbations in their coefficients [23]. To investigate this behavior, we conducted experiments on a perturbed version of the Wilkinson polynomial defined as $\prod_{i=1}^{n}(x+i) - 2^{23}x^{19}$. The accuracy results obtained using both the Durand-Kerner method with initial approximation $\lambda_{max}$ Bound and MATLAB's built-in *roots*() function were consistent with those reported by Wilkinson, as illustrated in Figure 5(a). However, the convergence behavior was significantly affected, with the algorithm reaching the maximum iteration threshold of 1000 without achieving full convergence, as shown in Figure 5(b).

Table 1. Accuracy of propose bound for polynomial Wilkinson

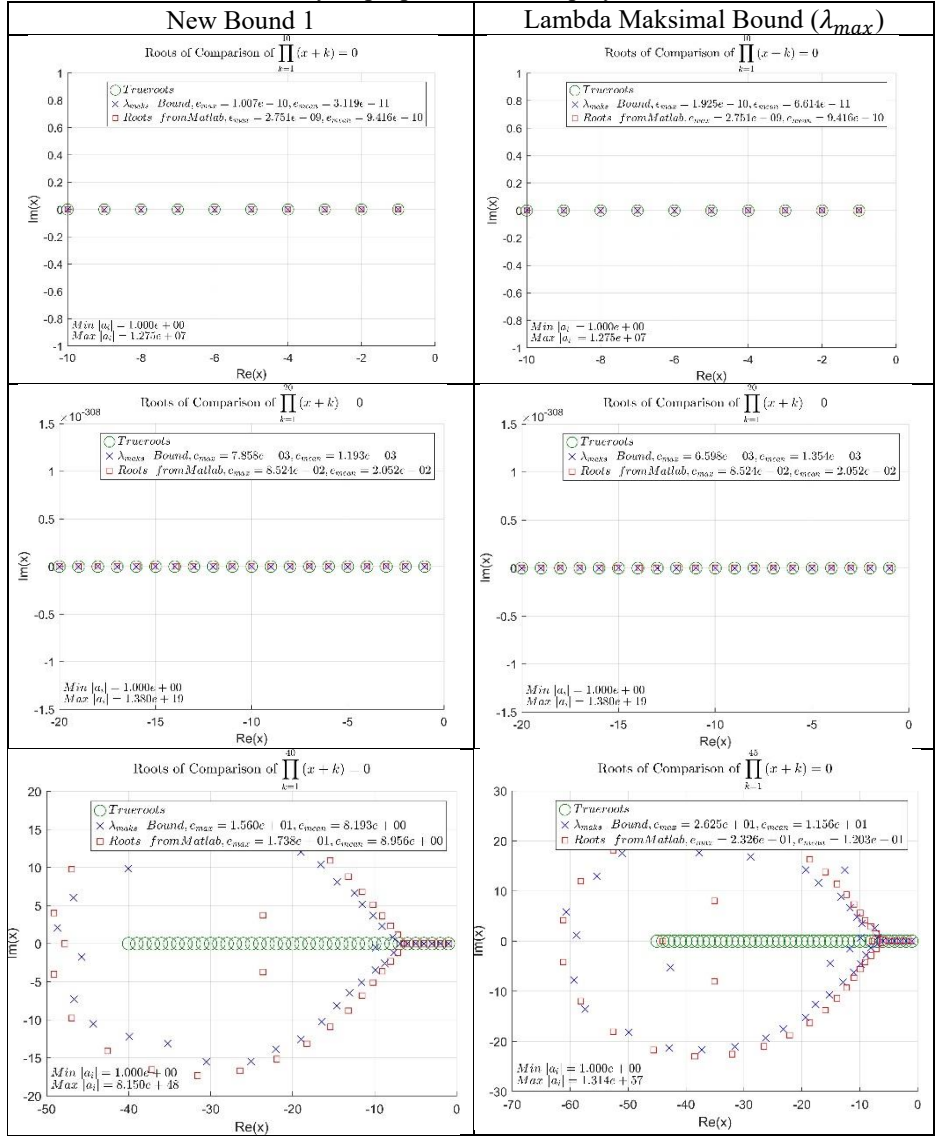

Tabel 2. Convergence of propose bound for polynomial Wilkinson

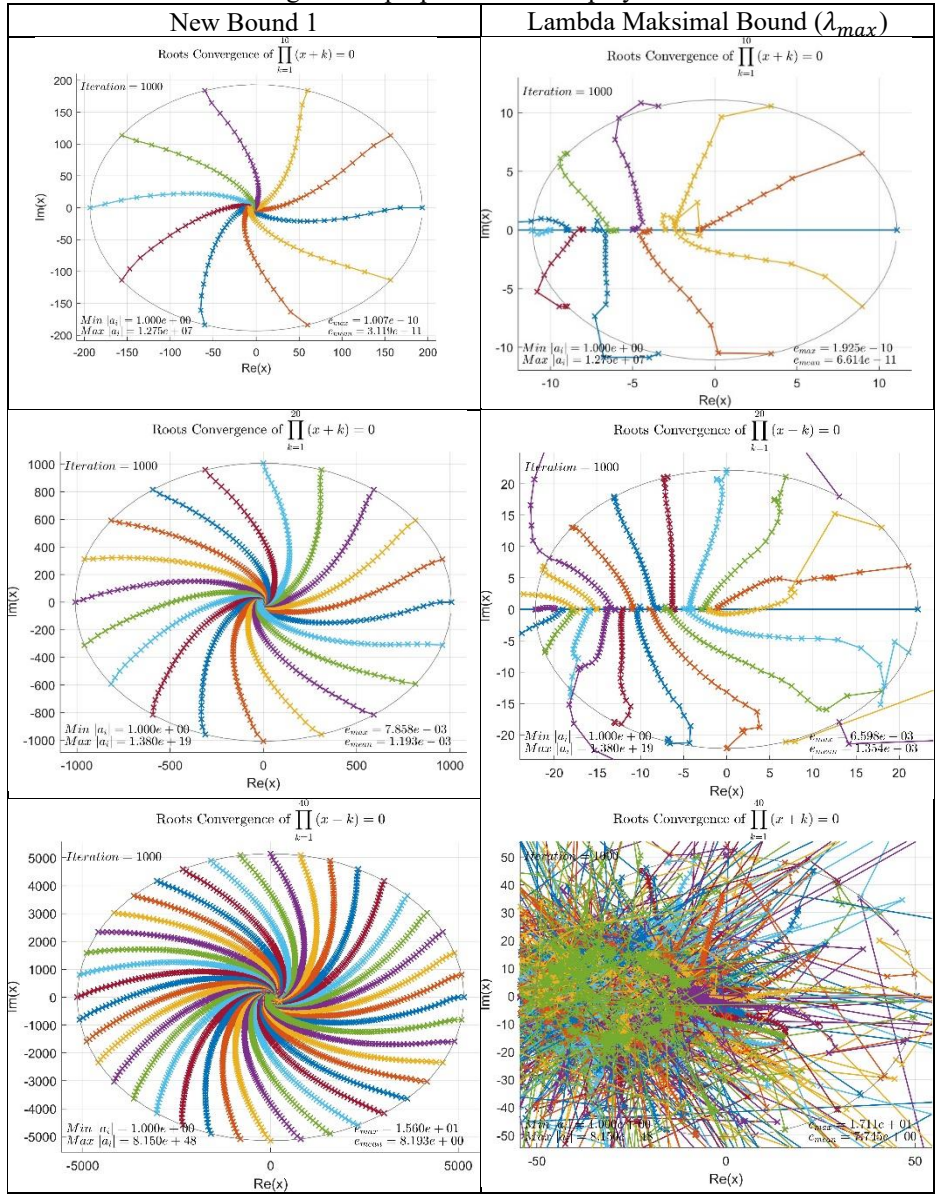

**Analysis on Clustered and Randomized Polynomials**

This subsection presents a detailed evaluation of the Durand-Kerner algorithm applied to polynomials characterized by clustered roots and randomly generated coefficients, two configurations known to exhibit numerical instability and convergence challenges. The analysis focuses on the algorithm's accuracy and convergence behavior when initialized using New Bound 1 and the Lambda Maximal Bound (λ_max), with comparative reference to MATLAB's built-in *roots*() function.

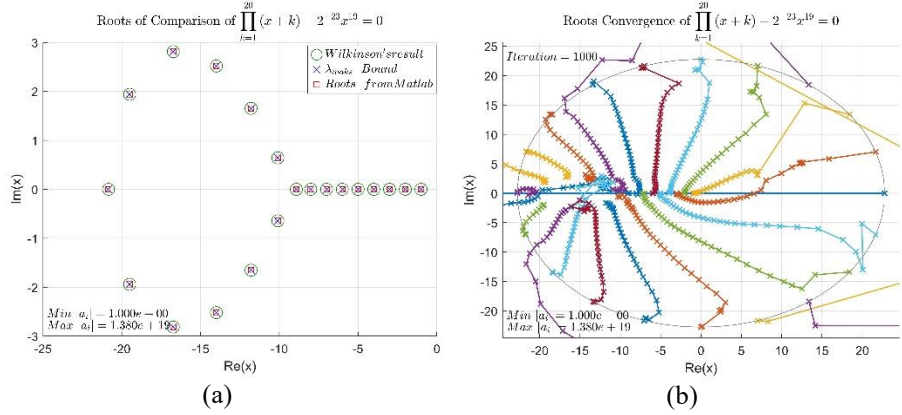

(a)  (b)

*Fig. 5. Accuracy and convergence of Wilkinson's polynomial with pertubation*

In the numerical experiments involving polynomials with clustered roots, we first employed the polynomial defined as $\prod_{i=1}^{n}(x + 1.000i)$. The accuracy of the Durand-Kerner algorithm using initial approximations derived from New Bound 1 and the Lambda Maximal Bound ($\lambda_{max}$) is illustrated in the plots presented in Table 3. For polynomials of degree up to $n = 30$, both bounds yielded high accuracy, with a mean error of $e_{mean} = 5.332 \times 10 - 1$, which is significantly better than the mean error produced by MATLAB's built-in *roots*() function, recorded at $6.146 \times 10 - 1$.

The convergence experiments were conducted using the clustered polynomial defined as $\prod_{i=1}^{n}(x + 1.000i)$. The results, illustrated in the plots presented in Table 4, evaluate the convergence behavior of the Durand-Kerner algorithm when initialized with New Bound 1 and the $\lambda_{max}$ Bound. For polynomial degrees up to $n = 30$, both bounds successfully guided the algorithm to converge to the true roots, comparable to the performance of MATLAB's built-in *roots*() function.

To investigate experiments on polynomials with random coefficients. We used random coefficient for $a_i, i = 2,3, …, n + 1$ in range of $[-15, 15]$. The accuracy and convergence behavior was significantly good for polynomial degree reaching $n = 140$, as shown in Figure 6.

To evaluate the performance of the Durand-Kerner algorithm on polynomials with randomly generated coefficients, we conducted experiments in which the coefficients $a_i$, for $i = 2,3, …, n + 1$, in the interval $[-15, 15]$. The results, as illustrated in Figure 6, demonstrate that both the accuracy and convergence behavior of the algorithm remain robust even for polynomials of high degree, up to $n = 140$.

Table 3. Accuracy of propose bound for polynomial with clustered roots

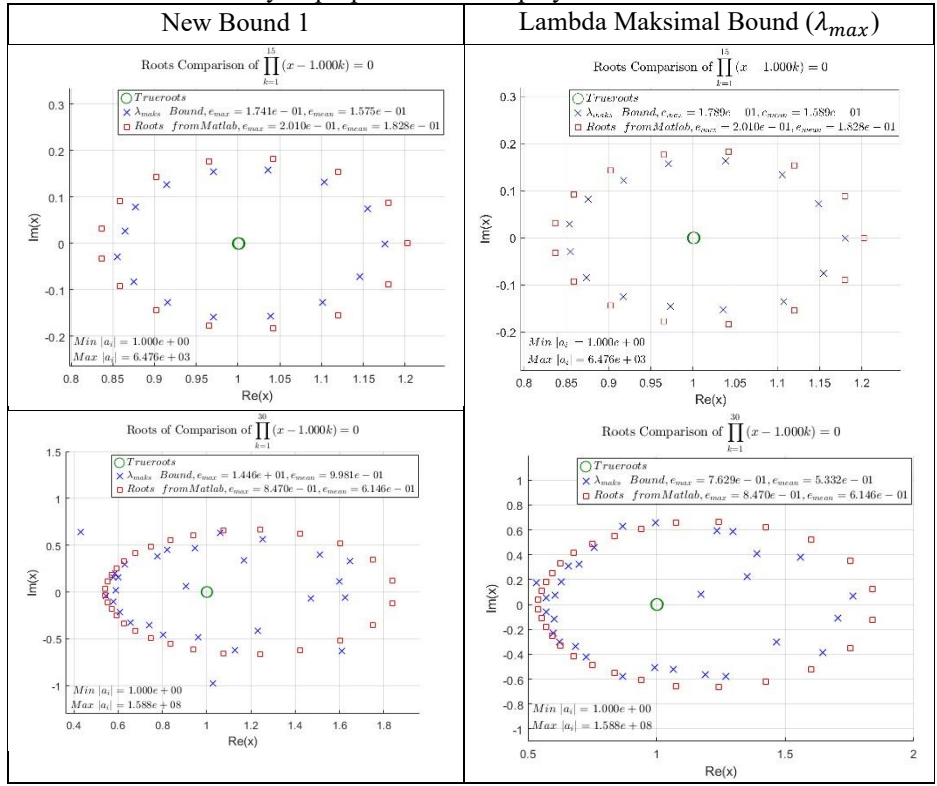

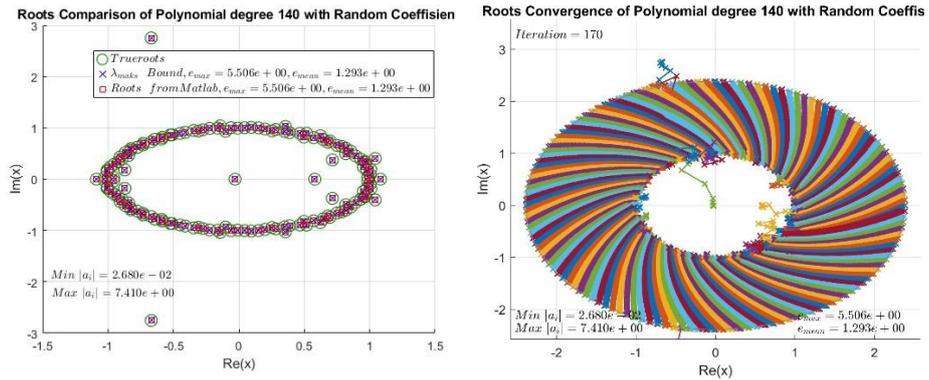

*Fig. 6. Accuracy and convergence of roots finding for polynomial ith random coefficients.*

Tabel 4. Convergence of propose bound for polynomial with clustered roots

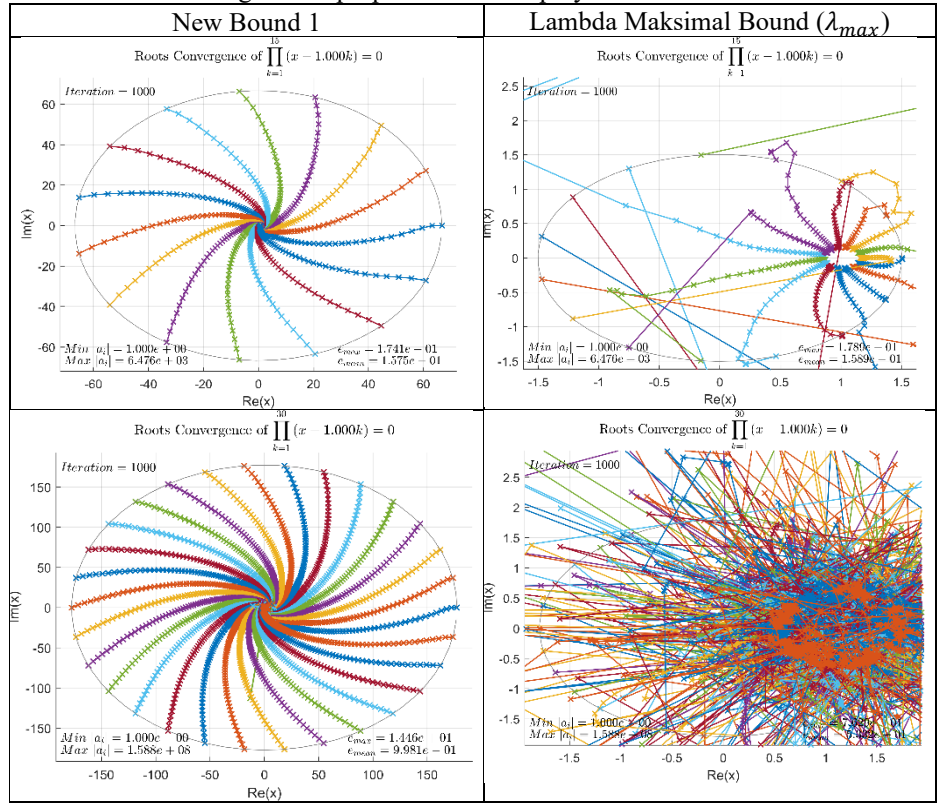

## Conclusion

This study has proposed and evaluated several initial approximations of the radius bounds for the Durand-Kerner algorithm in finding all the roots of a polynomial. Through the numerical experiments, the following conclusions can be drawn:
1. The coice of initial approximation selection significantly influences the convergence behavior. The choice of the initial radius directly affects the speed of convergence of the Durand-Kerner method.
2. New bound 1, while ensuring convergence, tends to produce excessively large radius. This can lead to inefficiencies in computation and more iterations.
3. The lambda maximal bound, derived from the dominant eigenvalue of the companion matrix, consistently guarantees that all roots lie within the complex circle and supports fast and stable convergence.

Overall, the proposed bounds enhance the robustness and efficiency of the Durand-Kerner method, particularly for polynomials with high degrees or challenging root structures.

## Future Work

Future research may focus on extending the proposed initial bound strategies to handle high-degree polynomials. Furthermore, integrating the initial approximation of the radius using machine learning techniques could further enhance the robustness and efficiency of the Durand-Kerner method.

Moreover, the Durand-Kerner algorithm is inherently parallelizable due to the independence of roots computation. Future work could explore the parallel algorithm based on multi-core CPUs or GPUs to reduce computation time, especially for high-degree polynomials. This parallelization could be combined with the proposed initial bounds to create efficient root-finding.

## Acknowledgment

We would like to thank the Department of Mathematics at Institut Teknologi Sepuluh Nopember (ITS), Surabaya, for providing the resources and assistance that allowed us to this research. Our appreciation also goes to colleagues who's insightful helped to make this research better.